\documentclass{hha}
\usepackage{amssymb,latexsym,amscd,verbatim}
\usepackage[all]{xy}

\newtheorem{hypothesis}[theorem]{Hypothesis}
\theoremstyle{remark}
\newtheorem{remarks}[theorem]{Remarks}
\newtheorem{question}[theorem]{Question}

\renewcommand{\epsilon}{\varepsilon}
\renewcommand{\phi}{\varphi}
\newcommand{\sA}{\mathcal{A}}
\newcommand{\sB}{\mathcal{B}}
\newcommand{\sC}{\mathcal{C}}
\newcommand{\sD}{\mathcal{D}}
\newcommand{\sE}{\mathcal{E}}

\newcommand{\sJ}{\mathcal{J}}
\newcommand{\bF}{\mathbf{F}}
\newcommand{\R}{\mathbf{R}}

\newcommand{\uI}{\underline{I}}

\newcommand{\by}[1]{\overset{#1}{\longrightarrow}}

\newcommand{\iso}{\by{\sim}}
\newcommand{\inj}{\hookrightarrow}

\newcommand{\bHom}{\operatorname{\mathbf{Hom}}}
\newcommand{\Ker}{\operatorname{Ker}}
\newcommand{\Var}{\operatorname{Var}}

\newcommand{\Sm}{\operatorname{Sm}}
\newcommand{\Sch}{\operatorname{Sch}}
\newcommand{\Spec}{\operatorname{Spec}}
\renewcommand{\lim}{\varprojlim}

\newcommand{\Coprod}{{\coprod}}
\newcommand{\proper}{{\operatorname{prop}}}
\newcommand{\proj}{{\operatorname{proj}}}

\newcommand{\op}{{\operatorname{op}}}
\newcommand{\qp}{{\operatorname{qp}}}
\newcommand{\sm}{{\operatorname{sm}}}

\newcommand{\cf}{{\it cf. }}
\newcommand{\rr}{\operatornamewithlimits{\rightrightarrows}\limits}

\newcounter{spec}
\newenvironment{thlist}{\begin{list}{\rm{(\roman{spec})}}%
{\usecounter{spec}\labelwidth=20pt\itemindent=0pt\labelsep=10pt}}%
{\end{list}}%

\numberwithin{equation}{section}

\newcommand{\resp}{{\it resp.} }
\newcommand{\car}{\operatorname{char}}

\begin{document}
\title{A few localisation theorems}
\author{Bruno Kahn}
\email{kahn@math.jussieu.fr}
\address{Institut de Math\'ematiques de Jussieu\\175--179 rue du
  Chevaleret\\75013 Paris\\France}
\author{R. Sujatha}
\email{sujatha@math.tifr.res.in}
\address{Tata Institute of Fundamental Research\\Homi Bhabha
  Road\\Mumbai 400 005\\India}
\begin{abstract} Given a functor $T:\sC \to \sD$ carrying a class of morphisms $S
\subset \sC$ into a class $S' \subset \sD$, we give sufficient conditions in order
that $T$ induces an equivalence on the localised categories. These conditions are
in the spirit of Quillen's theorem A. We give some applications in algebaic and
birational geometry. 
\end{abstract}
\maketitle

\section*{Introduction}
Let $T:\sC\to \sD$ be a functor and $S\subset \sC$, $S'\subset \sD$ two 
classes of morphisms containing identities and stable under
composition, such that $T(S)\subseteq S'$. This induces the situation
\begin{equation}\label{eq0}
\begin{CD}
\sC@>T>> \sD\\
@V{P}VV @V{Q}VV\\
S^{-1}\sC@>\bar T>> {S'}^{-1}\sD
\end{CD}
\end{equation}
where $P$ and $Q$ are localisation functors. In this note,
we offer an answer to the following question.

\begin{question} Give sufficient conditions for $\bar T$
to be an equivalence of categories.
\end{question}

This answer, Theorem \ref{t0}, is in the
spirit  of Quillen's  theorem A
\cite[th. A]{quillen2} that we recall for
motivation: in the above situation, forgetting $S$
and $S'$, if for all $d\in \sD$ the category $d\backslash T$ (see \S
\ref{comma})  is
 $\infty$-connected, then $T$ is a weak equivalence. 

\subsubsection*{Background} In \cite[Th. 3.8]{ks1}, we proved that $\bar
T$ is an equivalence of categories when $\sD$ is the category of smooth
varieties over a field of characteristic $0$, $\sC$ is its full
subcategory consisting of smooth projective varieties, and we take for $S$
and
$S'$ either birational morphisms or ``stable birational morphisms" (i.e.
dominant morphisms such that the corresponding function field extension
is purely transcendental). When we started revising \cite{ks1}, it
turned out that we needed similar localisation theorems in other
situations. At this stage
it was becoming desirable to understand these
localisation theorems more abstractly, and indeed we got two
non-overlapping, technical (and very ugly) statements. 

The first author then discussed these results with Georges Maltsiniotis,
and they arrived at Corollary \ref{c2} below. Using Proposition
\ref{p1} below, one can easily see 
that the hypotheses of Corollary \ref{c2} are verified in the case of
Theorem 3.8 of 
\cite{ks1}. However, they are not verified in some of the other
geometric situations mentioned above.  

``Catching" the latter situations led to Theorem \ref{t0}. Thus we had
two sets of abstract hypotheses implying that $\bar T$ is an
equivalence of categories: 
\begin{itemize}
\item hypotheses (0), (1) and (2) of Theorem \ref{t0}.
\item hypotheses (0) and (1') of Corollary \ref{c2};
\end{itemize}

To crown all, Maltsiniotis gave us an argument showing that
$(0) + (1')\Rightarrow (0)+(1)+(2)$: this is the content of Theorem
\ref{t1} a) and the proof we give is essentially his. 

In the same period, Jo\"el Riou proved a localisation theorem
of a similar nature (Theorem \ref{riou}). It turns out that Hypotheses
(0), (1) and (2) are implied by Riou's hypotheses (and actually by
less): see Theorem \ref{p3}. 

After stating and proving the main theorem, Theorem \ref{t0}, we prove
a ``relativisation" theorem, Theorem \ref{t1} which leads to Corollary
\ref{c2} mentioned above. We then give a number of conditions which
imply the hypotheses of Theorem \ref{t0} in \S \ref{5}. In \S \ref{6}
we show that the fact that $\bar T$ is an equivalence of categories in
\eqref{eq0} is stable under adjoining products and coproducts. We then
give some algebro-geometric
      applications (hyperenvelopes, cubical hyperresolutions...) in \S
      \ref{7}, and finally, in \S \ref{8}, the birational applications
      we alluded to: those will be used to simplify the exposition of
      the revision \cite{ks} of \cite{ks1}. 

Even though Maltsiniotis did not wish to appear as a coauthor of this
note,  we want to stress his essential contributions in bringing the
results here to their present form. Let
us also mention that Hypotheses (0), (1) and (2) imply much more than Theorem
\ref{t0}: they actually yield the existence of an ``absolute" derived
functor (in the sense of Quillen \cite[\S 4.1, Def. 1]{quillen1})
associated to any functor $F:\sD\to\sE$ such that
$FT(S)$ is invertible. This will be developed in a forthcoming work of
Maltsiniotis and the first author, where a different proof of Theorem
\ref{t0} will be given \cite{KM}; see already \S \ref{kan} here for a
weaker result. In \cite{KM}, we also hope to lift Theorem \ref{t0} to
the ``Dwyer-Kan localisation" \cite{dk} by suitably reinforcing its
hypotheses. Finally, we wish to thank the referee for a very helpful comment
regarding Theorem \ref{p3} (see Lemma \ref{inftycon}).

\section{Notation} 

\subsection{Comma categories}\label{comma} Recall \cite[ch. II,
\S6]{mcl} that to a diagram of categories and functors
\[\begin{CD}
&&\sB\\
&&@VV{G}V\\
\sA@>F>> \sC
\end{CD}\]
one associates a category $F\downarrow G$, the (ordered) ``2-fibred
product" of $F$ and $G$: 
\[\begin{CD}
F\downarrow G@>F'>>\sB\\
@VV{G'}V@VV{G}V\\
\sA@>F>> \sC.
\end{CD}\]

An object of $F\downarrow G$ is a triple $(a,b,f)$ where $a\in \sA$,
$b\in \sB$ and $f$ is a morphism from $F(a)$ to $G(b)$. A morphism from
$(a,b,f)$ to $(a',b',f')$ is a pair of morphisms $\phi:a\to a'$,
$\psi:b\to b'$ such that the diagram
\[\begin{CD}
F(a)@>f>> G(b)\\
@VV{F(\phi)}V @VV{G(\psi)}V\\
F(a')@>f'>> G(b')
\end{CD}\]
commutes. Composition of morphisms is defined in the obvious way.

This notation is subject to the following abbreviations:
\begin{itemize}
\item $G=Id_\sB$: $F\downarrow G=F\downarrow \sB$.
\item Dually, $F=Id_\sA$: $F\downarrow G=\sA\downarrow G$.
\item If $\sB$ is the point category and $G$  has image $c$:
$F\downarrow G=F\downarrow c=F/c=\sA/c$ (the latter notation being
  used only when there is no possible ambiguity). 
\item Dually, if $\sA$ is the point category and $F$ has image
$c$: $F\downarrow G=c\downarrow G=c\backslash G=c\backslash \sB$.
\end{itemize}

The category $F\downarrow G$ should not be confused with its
full subcategory $F\times_\sC G$ or $\sA\times_\sC\sB$ ($1$-fibred
product), consisting of those triples $(a,b,f)$ such that $f$ is an
identity.

\subsection{Path groupoid} For any category $\sE$, one
denotes by $\Pi_1(\sE)$ the category obtained by inverting all arrows
of $\sE$: this is the \emph{path groupoid} of $\sE$. 

\subsection{Connectedness} A category $\sE$ is \emph{$n$-connected} if
(the geometric realisation  of) its nerve is $n$-connected;
$-1$-connected is synonymous to ``non-empty". For $n\le 1$, $\sE$ is
$n$-connected if and only if
$\Pi_1(\sE)$ is $n$-connected. Thus, $0$-connected means that $\sE$ is
nonempty and any two of its objects may be connected by a zig-zag of
arrows (possibly not all pointing in the same direction) and
$1$-connected means that $\Pi_1(\sE)$ is equivalent to the point (category
with one object and one morphism). 

If $\sE$ is
$n$-connected for any
$n$, we say that it is \emph{$\infty$-connected} (this notion is
apparently weaker than ``contractible").

\subsection{Cofinal functors}\label{fin} According to \cite[ch. IX, \S 3
p. 217]{mcl}, a functor $L:\sJ'\to \sJ$ is called \emph{cofinal} if, for
all $j\in \sJ$, the category 
$L\downarrow j=L/j$ is $0$-connected. 

\begin{comment}
One has the following theorem, dual to
\cite[loc. cit., th. 1 p. 217 and ex. 5 p. 218]{mcl}:

\begin{theorem} \label{tfin} Let $L:\sJ'\to \sJ$ be a functor between
  small categories. The following conditions are equivalent: 
\begin{thlist}
\item $L$ is cofinal.
\item For any functor $F:\sJ\to X$ such that $\lim_{\sJ'} FL$ exists,
  the limit $\lim_\sJ F$ exists and the canonical morphism 
\[\lim_{\sJ'} FL\to \lim_\sJ F\]
is an isomorphism.
\item Same statement by limiting to $X=Set$ and $F$ corepresentable.
\end{thlist}
\end{theorem}

We shall use this theorem in the proof of Theorem \ref{t1}.
\end{comment}

\section{The main localisation theorem}

\subsection{The categories $I_d$}\label{Id} With the notation of the
introduction, consider $S$ and $S'$ as subcategories of $\sC$ and $\sD$
with the same objects, and let $T_S:S\to S'$ be the functor induced by
$T$. Set, for all
$d\in\sD$,
\[I_d=d\downarrow T_S=d\backslash S\]
cf. \ref{comma}. Thus:
\begin{itemize}
\item An object of $I_d$ is a pair $(c,s)$ where $c\in \sC$ and $s:d\to
T(c)$ belongs to $S'$. We summarise this with the notation $d\by{s}
T(c)$, or sometimes $s$, or even $c$ if this does not cause any confusion.
\item If $d\by{s} T(c)$, $d\by{s'} T(c')$ are two objects of $I_d$, a
morphism from the first to the second is a morphism $\sigma:c\to c'$
belonging to $S$ and such that the diagram 
\[\xymatrix{
&T(c)\ar[dd]^{T(\sigma)}\\
d\ar[ur]^s\ar[dr]^{s'}\\
&T(c')
}\]
commutes, composition of morphisms being the obvious one.
\end{itemize}

\subsection{Categories of diagrams}\label{diagrams} Let $E$ be a  small
category. In the category $\sC^E=\bHom(E,\sC)$, one
may consider the following class of morphisms $S(E)$: if $c,c'\in
\sC^E$, a morphism
$s:c\to c'$ belongs to $S(E)$ if and only if, for all $e\in E$,
$s(e):c(e)\to c'(e)$ belongs to $S$. One defines similarly ${S'}(E)$,
a class of morphisms in $\sD^E$. This gives a meaning to the notation
$I_d$ for $d\in \sD^E$.

We shall be interested in the case where $E=\Delta^n$, corresponding to the
totally ordered set $\{0,\dots, n\}$: so, $\sC^{\Delta^n}$
can be identified with the category of sequences of $n$ composable arrows
$(f_n,\dots, f_1)$ of $\sC$. For $n=0$, this is just the
category $\sC$.

\subsection{Statement of the theorem} With notation as in \S\S
\ref{Id} and \ref{diagrams}, it 
is the following:

\begin{theorem}[Simplicial theorem]\label{t0} Suppose the following
  assumptions verified: 
\begin{itemize}
\item[\rm (0)] For all $d\in \sD$, $I_d$ is $1$-connected.
\item[\rm (1)] For all $f\in \sD^{\Delta^1}$, $I_f$ is $0$-connected.
\item[\rm (2)] For all $(f_2,f_1)\in \sD^{\Delta^2}$, $I_{(f_2,f_1)}$
  is $-1$-connected. 
\end{itemize}
Then $\bar T$ is an equivalence of categories.
\end{theorem}

\subsection{Preparatory lemmas}\label{2.4} Before proving theorem \ref{t0}, we
shall establish a few lemmas. The first is trivial:

\begin{lemma}\label{l00} For all $d\in \sD$, the composite functor
\[I_d\to \sC\by{P} S^{-1}\sC,\]
where the first functor sends $d\by{s} T(c)$ to $c$, inverts all
arrows of $I_d$, hence induces a functor 
\[\bF(d):\Pi_1(I_d)\to S^{-1}\sC.\]
\end{lemma}

For $d\in \sD$ and for $c,c'\in \Pi_1(I_d)$, denote by $\gamma_{c,c'}$
the unique morphism from $c$ to $c'$, as well as its image in
$Ar(S^{-1}\sC)$ by the functor $\bF(d)$. Let $f:d_0\to d_1$ be a 
morphism of $\sD$. For $(c_1,c_0,g)\in Ob(I_{d_1})\times
Ob(I_{d_0})\times Ob(I_f)$, set
\[\phi_f(c_1,c_0,g)= \gamma_{c_1,rg}^{-1}\circ
g\circ\gamma_{c_0,dg}\in S^{-1}\sC(c_0,c_1)\] 
where $dg,rg$ denote respectively the domain and the range of $g$. If
$g,g'\in I_f$, a morphism $g\to g'$ yields a
commutative diagram
\[\begin{CD}
dg@>g>> rg\\
@V{\sigma}VV @V{\tau}VV\\
dg'@>g'>> rg'
\end{CD}\]
with $\sigma\in Ar(I_{d_0})$, $\tau\in Ar(I_{d_1})$. We then have
\begin{multline*}
\phi_f(c_1,c_0,g')=\gamma_{c_1,rg'}^{-1}\circ
g'\circ\gamma_{c_0,dg'}=\gamma_{c_1,rg'}^{-1}\circ\tau
\circ g\circ \sigma^{-1}\circ\gamma_{c_0,dg'}\\
= \gamma_{c_1,rg}^{-1}\circ g\circ\gamma_{c_0,dg}
=\phi_f(c_1,c_0,g)
\end{multline*}
in view of the fact that $\sigma=\gamma_{dg,dg'}$ and
$\tau=\gamma_{rg,rg'}$ in $S^{-1}\sC$. 

Since $I_f$ is $0$-connected, one deduces a canonical map
\[\phi_f:Ob(I_{d_1})\times Ob(I_{d_0})\to Ar(S^{-1}\sC)\]
such that $d\phi_f(c_1,c_0)=c_0$ and $r\phi_f(c_1,c_0)=c_1$. Observe the
formula
\begin{equation}\label{e11}
\phi_f(c'_1,c'_0)=\gamma_{c'_1,c_1}^{-1}\phi_f(c_1,c_0)\gamma_{c'_0,c_0}.
\end{equation}

In other words, $\phi_f$ defines a \emph{functor}
$\Pi_1(I_{d_0})\times \Pi_1(I_{d_1})\to (S^{-1}\sC)^{\Delta^1}$ lifting
the functors $\bF(d_0)$ and $\bF(d_1)$ via the commutative diagram
\[\xymatrix{
\Pi_1(I_{d_0})\times \Pi_1(I_{d_1})\ar[r]^{\phi_f}\ar[dr]_{\bF(d_0)\times
\bF(d_1)}& (S^{-1}\sC)^{\Delta^1}\ar[d]^{(d,r)}\\
&S^{-1}\sC\times S^{-1}\sC.
}\]

\begin{lemma}\label{l01}
a) If $f=1_d$ for some $d\in \sD$, then $\phi_f(c,c)=1_c$ for all $c\in
Ob(I_d)$.\\
b) If $f_1:d_0\to d_1$ and $f_2:d_1\to d_2$, then
\[\phi_{f_2f_1}(c_2,c_0)=\phi_{f_2}(c_2,c_1)\phi_{f_1}(c_1,c_0)\]
for all $(c_0,c_1,c_2)\in Ob(I_{d_0})\times Ob(I_{d_1})\times
Ob(I_{d_2})$.\\ 
c) If $f\in S'$, $\phi_f(c_1,c_0)$ is invertible in
$S^{-1}\sC$ for all $(c_0,c_1)\in Ob(I_{d_0})\times Ob(I_{d_1})$.
\end{lemma}

\begin{proof} a) is obvious. To prove b), let us use hypothesis
(2) to find $g_1:c_0\to c_1$ and $g_2:c_1\to c_2$ respectively in
$I_{f_1}$ and $I_{f_2}$. Then $\phi_{f_1}(c_1,c_0)=g_1$,
$\phi_{f_2}(c_2,c_1)=g_2$ and $\phi_{f_2f_1}(c_2,c_0)=g_2g_1$. Hence b) is
true for this particular choice of $(c_0,c_1,c_2)$, and one deduces from
\eqref{e11} that it remains true for all other choices.

Let us prove c). Choose a commutative diagram ($-1$-connectedness of
$I_f$)
\[\begin{CD}
d_0@>s_0>> T(c'_0)\\
@V{f}VV @V{T(g)}VV\\
d_1@>s_1>> T(c'_1)
\end{CD}\]
where $s_0,s_1\in S'$. Since $S'$ is stable under composition, this
diagram shows (using $s_1f$) that $g$ defines an object of
$I_{1_{d_0}}$; moreover, 
one obviously has $\phi_{1_{d_0}}(c'_1,c'_0)=g$. From a) and
\eqref{e11} (applied with $c_0=c_1$), one deduces that $g$ is invertible.
On the other hand, one also has
$g=\phi_f(c'_1,c'_0)$; reapplying \eqref{e11}, we get the desired
conclusion.
\end{proof}

\subsection{Proof of Theorem \ref{t0}}\label{2.5} We start by
defining a functor
\[F:\sD\to S^{-1}\sC\]
as follows: for all $d\in Ob(\sD)$, choose an object
$d\by{s_d}T(c_d)$ of $I_d$. Set
\begin{align*}
F(d)&=c_d\\
F(f)&=\phi_f(c_{d_1}, c_{d_0})\text{ for } f:d_0\to d_1.
\end{align*}

Lemma \ref{l01} shows that $F$ is indeed a functor, and that it inverts
the arrows of $S'$; hence it induces a functor
\[\bar F:{S'}^{-1}\sD\to S^{-1}\sC.\]

For $c\in Ob(S^{-1}\sC)$, one has an isomorphism
\[\gamma_{c,c_{T(c)}}:\bar F\bar T(c)\iso c.\]

Formula \eqref{e11} shows that it is natural in $c$: one checks it
first for the morphisms of $\sC$, then naturality passes 
automatically to $S^{-1}\sC$. On the other hand, for $d\in
Ob({S'}^{-1}\sD)$, one has an isomorphism
\[s_d:d\iso\bar T\bar F(d).\]

The definitions of $\phi_f$ and \eqref{e11} show again that this
isomorphism is natural in $d$ (same method).

It follows that $\bar F$ is quasi-inverse to $\bar T$.  \qed

\section{Towards Kan extensions}\label{kan}

\enlargethispage*{20pt}

Keep the setting of \eqref{eq0} and hypotheses of Theorem \ref{t0}, and let
$F:\sD\to\sE$ be another functor. We assume:

\begin{hypothesis}\label{h1} The functor $FT$ inverts $S$, i.e., there exists
  a functor $G:S^{-1}\sC\to \sE$ and a natural 
isomorphism
\[FT \simeq GP.\]
\end{hypothesis}

Under this hypothesis, let us define
\[RF:=G\bar T^{-1}:{S'}^{-1}\sD\to \sE\]
where $\bar T^{-1}$ is a chosen quasi-inverse to $\bar T$.

We construct a natural transformation $\eta:F\Rightarrow RF\circ Q$ as follows:

Let $d\in \sD$ and $d\by{s} T(c)\in I_d$. Then $s$ defines
\[\begin{CD}
F(d)@>{F(s)}>>FT(c)\simeq GP(c)@>{G\bar T^{-1}Q(s)^{-1}}>> G\bar
T^{-1} Q(d)= RF\circ Q(d).
\end{CD}\]

Since $I_d$ is $0$-connected, this morphism $\eta_d$ does not depend
on the choice 
of $s$. Then, the $-1$-connectedness of the categories $I_f$ shows
that $\eta$ is 
indeed a natural transformation.

It will be shown in \cite{KM} that $(RF,\eta)$ is in fact a left Kan extension
\cite[ch. X, \S 3]{mcl} (= right total derived functor \`a la Quillen
\cite[\S 4.1, Def. 1]{quillen1}) of $F$ along $Q$, but this requires the
full force of the hypotheses of Theorem \ref{t0}.

\section{A relativisation theorem}

\subsection{Two lemmas on comma categories}

\begin{lemma}[``theorem a"]\label{lcomma1} Let $F:\sA\to \sB$ be a final
functor (\S \ref{fin}). Then $F$ induces a bijection on
the sets of 
connected components. In particular, $\sA$ is
$0$-connected if and only if $\sB$ is $0$-connected.
\end{lemma}

\begin{proof} (See also \cite[Ex. 1.1.32]{maltsin}.) Surjectivity is 
obvious. For injectivity, let $a_0,a_1\in \sA$ be such that
$F(a_0)$ and $F(a_1)$ are connected. By the surjectivity of $F$, any
vertex of a 
chain linking them is of the form $F(a)$. To prove
that $a_0$ and $a_1$ are connected, one can therefore reduce to the case where
$F(a_0)$ and $F(a_1)$ are directly connected, say by a 
morphism $f:F(a_0)\to F(a_1)$. But the two objects
\[F(a_0)\by{f}F(a_1),\qquad F(a_1)\by{=}F(a_1)
\]
of $F/F(a_1)$ are connected by assumption, which implies that $a_0$ and $a_1$
are connected in $\sA$.
\end{proof}

\begin{lemma}\label{lcomma2} Let
\[\begin{CD}
F\downarrow G@>F'>>\sB\\
@VV{G'}V@VV{G}V\\
\sA@>F>> \sC.
\end{CD}\]
be a ``$2$-cartesian square" of categories.\\
a) For all $b\in \sB$, the functor
\begin{align*}
G_*:F'/b\to F/G(b)\\
\left[
\begin{CD}
F(a)@>f>> G(b')\\
&&@V{G(\phi)}VV\\
&& G(b)
\end{CD}
\right]
&\mapsto\quad [F(a)\by{G(\phi)f} G(b)]
\end{align*}
has a right adjoint/right inverse $G^!$ given by
\[G^!([F(a)\by{f} G(b)])=
\left[
\begin{CD}
F(a)@>f>> G(b)\\
&&@V{G(1_b)}VV\\
&& G(b).
\end{CD}
\right]\]
In particular, $G_*$ is a weak equivalence.\\
b) Suppose that $F/c$ is nonempty for all $c\in \sC$. Then $F'$ is
surjective.\\
c) Suppose moreover that $F$ is cofinal. Then $F'$ induces a bijection on
connected components.
\end{lemma}

\begin{proof} a) is checked immediately; the fact that $G_*$ is a weak
equivalence then follows from \cite[p. 92, cor. 1]{quillen2}. b) is
obvious. It remains to prove c): by a), the
categories $F'/b$ are $0$-connected. The conclusion then follows from
Lemma \ref{lcomma1} applied to $F'$.
\end{proof}

\subsection{The theorem} For all $d\in \sD$, let now $J_d:=d\backslash
\sD$.

\begin{theorem}\label{t1} a)  Suppose the following
conditions hold for all
$d\in Ob(\sD)$:
\begin{itemize}
\item[\rm (0)] $I_d$ is $1$-connected.
\item[\rm (1')] The obvious functor $\Phi_d:I_d\to J_d$ is cofinal (\S
\ref{fin}).
\end{itemize}
Then, for all $n\ge 0$ and all $(d_0\to\dots\to d_n)\in \sD^{\Delta^n}$,
the category $I_{(d_0\to\dots\to d_n)}$ is $0$-connected. \\  
b) Suppose that, for all $d\in \sD$ all $j\in J_d$, $I_d$ and
$I_d/j$ are  $\infty$-connected. Then, for all $n\ge 0$ and all 
$(d_0\to\dots\to d_n)\in \sD^{\Delta^n}$, the category
$I_{(d_0\to\dots\to d_n)}$
is $\infty$-connected.
\end{theorem}

\begin{proof}  a) One proceeds by induction on $n$, the case $n=0$
following from Hypothesis (0). Consider the obvious forgetful functors
\[I_{(d_0\to\dots\to d_n)}\by{u} I_{(d_1\to\dots\to d_n)}\by{v} I_{d_1}, \qquad
I_{(d_0\to\dots\to d_n)}\by{w} I_{d_0}.\]
 
One checks immediately that the diagram
\begin{equation}\label{eqcomma}
\begin{CD}
I_{(d_0\to\dots\to d_n)}@>u>> I_{(d_1\to\dots\to d_n)}\\
@V{w}VV @VV{f_1^*\circ\Phi_{d_1}\circ v}V\\
I_{d_0}@>\Phi_{d_0}>> J_{d_0}.
\end{CD}
\end{equation}
induces an isomorphism of categories
\[I_{(d_0\to\dots\to d_n)}=\Phi_{d_0}\downarrow(f_1^*\circ\Phi_{d_1}\circ v)\]
i.e. is $2$-cartesian. Here, $f_1:d_0\to d_1$. Hypothesis (1') then
implies that Lemma \ref{lcomma2} 
c) can be applied with $F=\Phi_{d_0}$. Therefore $u$ induces a
bijection on connected components, hence the conclusion.

b)  Let us use Diagram \eqref{eqcomma} again. It follows from Lemma
\ref{lcomma2} a) that $u/x$ is 
$\infty$-connected  for all $x\in I_{(d_1\to\dots\to d_n)}$. By Quillen's
theorem A \cite[th. A]{quillen2}, $u$ is a weak equivalence; by
induction on $n$,
$I_{(d_1\to\dots\to d_n)}$ is  $\infty$-connected, hence so is
$I_{(d_0\to\dots \to d_n)}$.
\end{proof}

\begin{corollary}[Normand theorem]\label{c2} Suppose the following
conditions \allowbreak verified for all
$d\in Ob(\sD)$:
\begin{itemize}
\item[\rm (0)] $I_d$ is $1$-connected.
\item[\rm (1')] The obvious functor $\Phi_d:I_d\to J_d$ is cofinal (\S
\ref{fin}).
\end{itemize}
Then $\bar T$ is an equivalence of categories.
\end{corollary}

\begin{proof} This follows from Theorem \ref{t1} a) and Theorem \ref{t0}.
\end{proof}

\begin{remark} There is an $n$-connected version of Quillen's
theorem A for any $n$ (cf. Maltsiniotis \cite[1.1.34]{maltsin},
Cisinski  \cite{cisinski}). Using it, 
one may replace $\infty$-connected by $n$-connected in the hypothesis and
conclusion of Theorem \ref{t1} b) (same proof).
\end{remark}

\section{Complements}\label{5}

\subsection{A relative version}

\begin{corollary}\label{c1}  Suppose that $T$ is fully
faithful.\\
a) If Conditions {\rm (0), (1), (2)} of Theorem \ref{t0} are
satisfied, they are also satisfied for all
$c\in\sC$ for the functor $c\backslash\sC\to
  T(c)\backslash\sD$ induced by $T$.\\
b) Same result with Conditions {\rm (0), (1')} of Theorem \ref{t1} a).\\
In particular, in case a) or b), the functor
\[S^{-1}(c\backslash\sC)\to {S'}^{-1}(T(c)\backslash\sD)\]
induced by $T$ is an equivalence of categories. 
\end{corollary}

\begin{proof} For $\delta=[T(c)\to d]\in
  T(c) \backslash\sD$, the full faithfulness of $T$ implies that the
forgetful functors
\[\delta\backslash(c\backslash\sC)\to d\backslash\sC,\qquad
  \delta\backslash(c\backslash S)\to d\backslash S\]
are isomorphisms of categories. Similarly when dealing with
categories of type $I_f$ and $I_{(f_2,f_1)}$.
\end{proof}

\subsection{Riou's theorem}
A statement similar to Corollary \ref{c2} was obtained independently by 
Jo\"el Riou: 

\begin{theorem}[Riou \protect{\cite[II.2.2]{riou}}]\label{riou} Suppose that
\begin{thlist}
\item $T$ is fully faithful; $S= S'\cap \sC$.
\item In $\sD$, push-outs of arrows of $S'$ exist and belong to $S'$.
\item If $s\in S'$ and the domain of $s$ is in $T(\sC)$, so is its range.
\item For any $d\in \sD$, $I_d\neq \emptyset$.
\end{thlist}
Then $\bar T$ is an equivalence of categories.
\end{theorem} 

(Riou's hypotheses are actually dual to these: we write them as above for
an easy comparison with the previous results. Also, Riou does not
assume that $S'$ is stable 
under composition.)

Riou's proof is in the style of that of Theorem \ref{t0}, but more
direct because 
push-outs immediately provide a functor. Actually, as we realised when
reading Gillet--Soul\'e \cite{gs}, his hypothesis (iii) is not 
necessary, as is shown by the following

\begin{theorem}\label{p3} Assume that the hypotheses {\rm (i), (ii), (iv)}
  of Theorem \ref{riou} 
are verified. Then:\\ 
a) For any finite partially ordered set $E$, these hypotheses are verified for
$T^E:\sC^E\to \sD^E$ and $S(E), S'(E)$ (cf. \S\ref{diagrams}).\\ 
b) For any $d\in \sD$, $I_d$ is $1$-connected (and even $\infty$-connected, see Lemma \ref{inftycon}).\\ 
c) In the situation of a), the hypotheses of Theorem \ref{t0} are verified;
in particular, $\bar T^E$ is an equivalence of categories.
\end{theorem}

\begin{proof} a) It suffices to prove (iv): for this, we argue by
  induction on $|E|$, the case $E=\{0\}$ being Hypothesis (iv). 

Suppose that $|E|>0$, and let $d_\bullet\in \sD^E$. Let $e\in E$ be a
maximal element, $E'=E-\{e\}$ and $d'_\bullet$ the restriction of
$d_\bullet$ to $\sD^{E'}$. By induction, pick an object
$d'_\bullet\by{s'_\bullet} T^{E'}(c'_\bullet)$ in $I_{d'_\bullet}$. 

Let $F$ be the set of those maximal elements of $E'$ which are
$<e$. If $F=\emptyset$, we just pick $d_e\by{s_e} T(c_e)$ in $I_{d_e}$
(by (iv)) and adjoin it to the previous object. If $F$ is not empty,
let $d'$ be the push-out of the maps $d_f\by{s_f} T(c_f)$ (for $f\in
F$) along the maps $d_f\to d_e$. By (ii), the map $d_e\to d'$ is in
$S'$. Pick $d'\to T(c_e)$ in $I_{d'}$ by (iv), and define $s_e$ as the
composition $d_e\to d'\to T(c_e)$. By (i), the compositions  
\[T(c_f)\to d'\to T(c_e)\]
define morphisms $\sigma_{f,e}:c_f\to c_e$ in $S$, and we are
done. (In picture: 
\[\xymatrix{
d_{f}\ar[r]^{s_{f}}\ar[d]_{f_{f,e}}& T(c_{f})\ar[d]\ar[dr]^{T(\sigma_{f,e})}\\
d_e\ar[r]\ar@/d3ex/[rr]^{s_e}& d'\ar[r]& T(c_e).)
}\]

b) Let $s,s'\in I_d$. Taking their push-out,  we get a new object
$d'\in \sD$; applying (iv) to 
$I_{d'}$, we then get a new object $s''\in I_d$ and maps $s\to s''$,
$s'\to s''$. In 
particular, $I_d$ is $0$-connected.

A similar argument shows that the first axiom of calculus of fractions holds in
$I_d$ (for the collection of all morphisms of $I_d$). Therefore, in
$\Pi_1(I_d)$, any morphism may be written as $u_2^{-1}u_1$ for
$u_1,u_2$ morphisms of $I_d$. To prove 
that $I_d$ is $1$-connected, it therefore suffices to show that, given
two morphisms $u_1,u_2\in 
I_d$ with the same domain and range, $u_1$ and $u_2$ become equal in
$\Pi_1(I_d)$.  

The following proof is inspired from reading
\cite[pp. 139---140]{gs}. Let $s:d\to T(c)$ and $s':d\to T(c')$ be the
domain and range of 
$u_1$ and
$u_2$. Consider the push-out diagrams
\[\begin{CD}
d@>s>> T(c)\\
@V{s}VV @V{a}VV\\
T(c)@>a>> d'\\
@V{T(u_1)}VV @V{g_1}VV \\
T(c')@>f>> d''
\end{CD}
\qquad
\begin{CD}
d@>s>> T(c)\\
@V{s}VV @V{a}VV\\
T(c)@>a>> d'\\
@V{T(u_2)}VV @V{g_2}VV \\
T(c')@>f>> d''.
\end{CD}
\]

Here $d''$ and $f$ are common to the two diagrams because $T(u_1)s =
T(u_2)s$. For the same reason, we have $g_1a=g_2a$ (vertically), hence
(in the lower squares) 
\[fT(u_1) = g_1a = g_2a = fT(u_2).\]
Choose $d''\by{s''} T(c'')$ in $I_{d''}$ by 1). Then $s''f= T(\sigma)$
for some $\sigma\in S$. Hence 
\[\sigma u_1 = \sigma u_2\]
and $u_1 = u_2$ in $\Pi_1(I_d)$, as requested.

c) For any $n\ge 0$, consider the ordered set $E\times \Delta^n$. Then
a) and b) show that, for any $d_\bullet\in \sD^{E\times \Delta^n} =
(\sD^E)^{\Delta^n}$, $I_{d_\bullet}$ is $1$-connected. In particular,
the hypotheses of Theorem \ref{t0} hold. 
\end{proof}

\begin{remark} It is not clear whether the conditions of Theorem \ref{riou} imply
Condition (1') of Theorem
\ref{t1} a).
\end{remark}

\begin{remark} We shall use Theorem \ref{p3} in the geometric applications.
\end{remark}

\begin{remark}
Even though the categories $I_d$ are 1-connected under the conditions of
Theorem
\ref{riou}, they are not filtering in general (for example, they are
not filtering 
in the geometric case considered by Riou). A natural 
question is whether they are $\infty$-connected. We would like to thank 
the referee for providing an affirmative answer and sketching an 
argument, which we reproduce in the lemma  
below.
\end{remark}

\begin{lemma}[Referee's lemma]\label{inftycon}
Let ($\sC,\sD,T,S,S')$ be as in the introduction. Suppose that, for any finite partially ordered set (poset) $E$ and for any $d_\bullet\in \sD^E$, the category $I_{d_\bullet}$ is $0$-connected. Then, for any finite poset $E$ and any $d_\bullet\in \sD^E$, $I_{d_\bullet}$ is $\infty$-connected.
\end{lemma}

\begin{proof}
We will use the following sufficient condition for a simplicial set to 
be $\infty$-connected: if $X$ is a simplicial set such that any map from 
the nerve of a finite partially ordered set (poset) is simplicially 
homotopic to a constant map, then $X$ is $\infty$-connected. (This can be proven, for example, using the fact that for any integer $k\ge 1$ the iterated subdivision $Sd^k(\partial\Delta^n)$ is the nerve of a finite poset and from the fact that Kan's $Ex^\infty(X)$ is a fibrant replacement of $X$.) 

From this,  one deduces that if $\sC$ is a small category such that for any finite 
poset $E$, the category of functors $\sC^E$ is 0-connected, then $\sC$ 
is $\infty$-connected (use the fact that the functor ``nerve" is fully faithful and commutes with finite products). Let now $d\in \sD$. We note that for a finite poset $E$, a functor 
$u$ from $E$ to $I_d$ is the same as a functor $v$ from $E$ to $S$ with 
a morphism of functors from $d_E$ to $T^E_S(v)$, where $d_E$ denotes the 
constant functor from $E$ to $S'$, with value $d$. In other words, we 
have an equivalence of categories $I_d^E \simeq I_{d_E}$. Hence by {\it 
a)} of Theorem \ref{p3}, we can apply {\it b)} to $T^E$ and conclude 
that $I_d^E$ is 0-connected, which proves that $I_d$ is 
$\infty$-connected. 

We may then apply this conclusion to the collection 
\[(\sC^E,\sD^E,T^E,S(E),S'(E))\] 
for any finite poset $E$, since for another finite poset $F$ we have $(\sC^E)^F\simeq \sC^{E\times F}$, etc.
\end{proof}

\subsection{Weakening the hypotheses} This subsection has grown out of
exchanges with 
Maltsiniotis.

For $d\in \sD$, let us write $J_d=d\backslash T$ as in
Theorem \ref{t1}. The projections $J_d\to \{d\}\subset\sD$ define a
fibred category $J$ over 
$\sD$. Similarly, the $I_d$ define a fibred category $I$ over $S'$ (viewed as
before as a category). 

Now replace $S$ and $S'$ by their strong saturations $\langle
S\rangle$ and $\langle 
S'\rangle$. (Recall that the strong saturation $\langle S\rangle$ of
$S$ is the collection, 
containing $S$, of morphisms $u\in \sC$ such that $u$ becomes
invertible in $S^{-1}\sC$.) We 
have similarly a fibred category
$\langle I\rangle$ over
$\langle S'\rangle$. For any $d$, we have obvious inclusions
\[I_d\subseteq \langle I\rangle_d\subseteq J_d.\]

We are interested in a collection of subcategories $I'_d$ of $\langle
I\rangle_d$ which form a 
fibred category over $S'$. Concretely, this means that, for any
$s:d\to d'$ in $S'$, the 
pull-back functor
\[s^*:\langle I\rangle_{d'}\to \langle I\rangle_d\]
sends $I'_{d'}$ into $I'_d$.

\begin{definition} A fibred category $I'\to S'$ as above is called a
  \emph{weak replacement} of $I$. 
\end{definition}

If $E$ is a small category, we have the fibred category $I(E)$ over
$S'(E)$ and we define a 
weak replacement of $I(E)$ similarly: namely, a collection of subcategories
$I'(E)_{\underline{d}}$ of $\langle I(E)\rangle_{\underline{d}}$
respected by pull-backs under 
morphisms of $S(E)$.

\begin{theorem}[Variant of Theorem \protect{\ref{t0}}] \label{t1v} 
Suppose given, for $n=0,1,2$, a weak replacement $I'(\Delta^n)$ of
$I(\Delta^n)$. Suppose 
moreover that 
\begin{itemize}
\item for any $f:d_0\to d_1$, the face functors $J_f\to J_{d_0}$ and
$J_f\to J_{d_1}$ send $I'_f$ to $I'_{d_0}$ and $I'_{d_1}$.
\item For any $(f_2,f_1)$, the face functors $J_{f_2,f_1}\to J_{f_2}$,
$J_{f_2,f_1}\to J_{f_1}$ and $J_{f_2,f_1}\to J_{f_2f_1}$ send
$I'_{(f_2,f_1)}$ respectively to $I'_{f_2}$, $I'_{f_1}$ and $I'_{f_2f_1}$.
\item For any $d\in \sD$, $I'_{1_d}$ contains at least one object of the
form $[1_d\to T(1_c)]$.
\end{itemize}
(The last condition is verified for example if the degeneracy functor
$J_d\to J_{1_d}$ sends $I'_d$ to $I'_{1_d}$.)\\
Finally, suppose that the $I'_{\underline{d}}$ have the same
connectivity properties as in 
Theorem \ref{t0}. Then $\bar T$ is an equivalence of categories. 
\end{theorem}

\begin{proof} One checks by inspection that the proof of Theorem \ref{t0} goes
through with these data.
\end{proof}

It was Maltsiniotis' remark that Corollary \ref{c2} still holds with a
weak replacement 
of $I$ rather than $I$. Presumably, one can check that Theorem
\ref{t1} still holds with weak 
replacements of the $I(\Delta^n)$, provided they satisfy simplicial
compatibilities similar to 
those of Theorem
\ref{t1v}.

\subsection{Sufficient conditions for (0), (1), (1') and (2)}

\begin{proposition}\label{p1} a) The following conditions imply the conditions 
of Theorem
\ref{t1} b) (hence, a fortiori, conditions {\rm (0)} and {\rm (1')} of Theorem
\ref{t1} a)): for any
$d\in\sD$ and
$j\in J_d$
\begin{itemize}
\item[(a1)] $I_d$ is cofiltering;
\item[(a2)] $I_d/j$ is (nonempty and) cofiltering.
\end{itemize}
b) The following conditions imply {\rm (a1)} and {\rm (a2)}:
\begin{itemize}
\item[(b1)] given $s\in S'$, $T(f)s=T(g)s\Rightarrow f=g$ ($f,g\in Ar(\sC)$);
\item[(b2)] $I_d$ is nonempty;
\item[(b3)] for any $(i,j)\in I_d\times J_d$, the $1$-fibred product
$I_d/i\times_{I_d}I_d/j$ is nonempty.
\end{itemize}
c) In b), conditions {\rm (b2)} and {\rm (b3)}
  are consequences of the following: finite products exist in $\sC$,
  $T$ commutes with them and, for any $d\in \sD$, there is a family of
objects
  $K_d\subset J_d$ such that 
\begin{itemize}
\item[(c1)] $K_d\ne\emptyset$; $I_d\subseteq K_d$; for any $k\in K_d$,
$I_d/k\ne\emptyset$.
\item[(c2)] If $k\in K_d$ and $j\in J_d$ then $j\times k\in K_d$. [Note
that the assumption on finite products implies that they exist in $J_d$
for any $d\in\sD$.]
\end{itemize}
\end{proposition}

\begin{proof} a) is ``well-known": see \cite[Prop. 2.4.9]{maltsin}. 

b) (b1)
  implies that $I_d$, hence also $I_d/j$, are ordered; (b2) and (b3) (the
latter
  applied with $j\in I_d$) then imply that $I_d$ is cofiltering and (b3)
  implies a fortiori that $I_d/j$ is nonempty for any $j\in J_d$; since
$I_d$ is cofiltering, $I_d/j$ is automatically cofiltering. 

c) Clearly (c1)
$\Rightarrow$ (b2). For (b3), let $(i,j)\in
I_d\times J_d$. By hypothesis,
$i\times j\in K_d$, hence $I_d/i\times j\neq\emptyset$ and there is an
$i'$ such that
$i'$ maps to $i\times j$, which exactly means that $i'\in
I_d/i\times_{I_d}I_d/j$.
\end{proof}

For the next proposition, we need to introduce a definition relative
to the pair $(\sD,S')$: 

\begin{definition}\label{d2} Given a diagram
\[\begin{CD}
d@>{s}>>d'\\
@V{f}VV\\
d_1
\end{CD}
\]
with $s\in S'$, we say that $s$ is \emph{in good position with respect
  to $f$} if the push-out 
\[\begin{CD}
d@>{s}>>d'\\
@V{f}VV@V{f'}VV\\
d_1@>{s_1}>>d'_1
\end{CD}\]
exists and $s_1\in S'$.
\end{definition}

\begin{proposition}\label{p2} Suppose that  the following
conditions  are verified:
\begin{itemize}
\item[(d1)] Morphisms of $S'$ are epimorphisms within $S'$.
\item[(d2)] If $f\in S'$ in Definition \ref{d2}, then any $s\in S'$ is
  in good position with respect to $f$. 
\item[(d3)] If $s\in S'$ is in good position with respect to $gf$,
  then it is in good position with respect to $f$. 
\item[(d4)] $T$ is fully faithful and $S=S'\cap \sC$.
\item[(d5)] For any $f:d\to d_1$ in $\sD$, there exists $s\in I_d$ in
  good position with respect to $f$. 
\end{itemize}
Then for all $m\ge 0$ and all $d_\bullet\in \sD^{\Delta^m}$,
$I_{d_\bullet}$ is ordered and filtering, hence $\infty$-connected. In
particular, the hypotheses of Theorem 
\ref{t0} are verified. 
\end{proposition}

\begin{proof} We first show that
$I_{d_\bullet}$ is nonempty. For $m=0$, this is (d5) applied to
  $f=1_{d_0}$. Suppose $m>0$: we argue by induction on $m$. Applying
  (d5) and (d3) to $f_m\circ\dots \circ f_1$, we find $s_0\in I_{d_0}$
  and a commutative (pushout) diagram 
\[\begin{CD}
d_0@>f_1>> d_1@>f_2>> \dots @>f_m>> d_m\\
@V{s_0}VV @V{s'_1}VV && @V{s'_m}VV\\
T(c_0)@>f'_1>> d'_1@>f'_2>> \dots @>f'_m>> d'_m
\end{CD}\]
with $s'_1,\dots, s'_m\in S'$. By induction, $I_{(d'_1\to\dots\to
  d'_m)}$ is nonempty, which shows using (d4) that $I_{d_\bullet}$ is
nonempty. (d1) then implies that it is ordered.  

Let us prove that they are filtering. Using (d2), we see that the
push-out $d_\bullet \to d'_\bullet$ of two objects
$d_\bullet\by{s_\bullet} T(c_\bullet)$, $d_\bullet\by{s'_\bullet}
T(c'_\bullet)$ of $I_{(f_m,\dots,f_1)}$ exists as a diagram in $\sD$;
using the nonemptiness of $I_{d'_\bullet}$, we conclude. 
\end{proof}

\begin{remark} This proposition (with its proof) may be seen as an easier
  variant of Theorem \ref{p3}. 
\end{remark}

\subsection{Another variant of Theorem \ref{t0}}

Keep notation as in Theorem \ref{t1}. As in \S \ref{Id}, let $S$ (\resp $S'$)
denote the subcategory of $\sC$ (\resp of $\sD$) with the same objects
but with only arrows in $S$ (\resp $S'$). Consider the category 
\[Id_S\downarrow T=\{(d,c,s)\mid d\in S',c\in S, s:d\to T(c)\}.\]

We have a projection functor 
\begin{gather*}
p_1:Id_S\downarrow T\to S'\\
(d,c,s)\mapsto d.
\end{gather*}

For $d\in S'$ we then define
\[\uI_d = p_1\downarrow d \]
so that objects of $\uI_d$ are diagrams
\begin{equation}\label{e2}
\begin{CD}
u@>s>> T(c)\\
@VjVV\\
d
\end{CD}
\end{equation}
with $s,j\in S'$, and morphisms are the obvious ones (in $S$).

We have the same definition for categories of diagrams. Then:

\begin{theorem}\label{tu0} Suppose the following assumptions verified:
\begin{itemize}
\item[\rm (\underline{0})] For all $d\in \sD$, $\uI_d$ is $1$-connected.
\item[\rm (\underline{1})] For all $f\in \sD^{\Delta^1}$, $\uI_f$ is
  $0$-connected. 
\item[\rm (\underline{2})] For all $(f_2,f_1)\in \sD^{\Delta^2}$,
$\uI_{(f_2,f_1)}$ is
$-1$-connected.
\end{itemize}
Suppose moreover that the following $2/3$ property holds:
\begin{itemize}
\item[\rm (*)] If $s\in S'$ and $st\in S'$, then $t\in S'$.
\end{itemize}
Then $\bar T$ is an equivalence of categories.
\end{theorem}

\begin{proof} One first mimics line by line the arguments of \S \ref{2.4}. The
only place where the added datum $j$ creates a difficulty is in the analogue of
Lemma \ref{l01} c). We then argue as follows: let $f:d_0\to d_1\in S'$. By the
$-1$-connectedness of $\uI_f$, we have a commutative diagram
\[\begin{CD}
d_0 @<j_0<< u_0@>s_0>> T(c'_0)\\
@V{f}VV @V{t}VV @V{T(g')}VV\\
d_1 @<j_1<< u_1@>s_1>> T(c'_1).
\end{CD}\]

Note that $j_1t=fj_0\in S'$, thus $t\in S'$ by (*), and therefore
$s_1t\in S'$. So 
we have another commutative diagram
\[\begin{CD}
d_0 @<j_0<< u_0@>s_0>> T(c'_0)\\
@V{1_{d_0}}VV @V{1_{u_0}}VV @V{T(g')}VV\\
d_0 @<j_0<< u_0@>s_1t>> T(c'_1)
\end{CD}\]
describing an object of $\uI_{1_{d_0}}$. From there, one proceeds as in the
proof of Lemma \ref{l01} c).

The analogue of
\S \ref{2.5} is now as follows: for each
$d\in
\sD$ one chooses an object
$(u_d,c_d,j_d,s_d)\in \uI_d$ and one defines a functor $F:\sD\to S^{-1}\sC$ by
$F(d)=c_d$, $F(f)=\phi_f(c_{d_1}, c_{d_0})$ as in \S \ref{2.5}. The natural
isomorphism $\bar F\bar T\Rightarrow Id_{S^{-1}\sC}$ is defined as in
\S \ref{2.5}; 
on the other hand, the isomorphism $Id_{{S'}^{-1}\sD}\Rightarrow \bar T\bar
F$ is defined on an object $d\in {S'}^{-1}\sD$ by $s_d j_d^{-1}$: it is easy to
check that it is natural.
\end{proof}

\section{Adding finite products or coproducts}\label{6}

In this section, we show that the property for $\bar T$ to be an
equivalence of categories in \eqref{eq0} is preserved by adjoining
finite products or coproducts. We shall only treat the case of
coproducts, since that of products is dual.

We shall say that a category $\sC$ \emph{has finite coproducts} (or
that $\sC$ is with finite coproducts) if all finite coproducts are
representable in $\sC$. This is the case if and only if
$\sC$ has an final object (empty coproduct) and the coproduct
of any two objects exists in $\sC$.

\begin{proposition}\label{p6.1}
Let $\sC$ be a category. There exists a category $\sC^\Coprod$
with finite coproducts and a
functor $I:\sC\to \sC^\Coprod$ with the following $2$-universal
property: any functor $F:\sC\to \sE$ where $\sE$ has finite coproducts
extends through $I$, uniquely up to natural isomorphism, to a functor
$F^\Coprod:\sC^\Coprod\to \sE$ which commutes with finite
coproducts. We call $\sC^\Coprod$ the \emph{finite coproduct
  envelope} of $\sC$.
\end{proposition}

\begin{proof} We shall only give a construction of $\sC^\Coprod$:
  objects are families $(C_i)_{i\in I}$ where $I$ is a finite set and
  $C_i\in \sC$ for all $i\in I$. A morphism $\phi:(C_i)_{i\in I}\to
  (D_j)_{j\in J}$ is given by a map $f:I\to J$ and, for all $i\in I$,
  a morphism $\phi_i:C_i\to D_{f(i)}$. Composition is defined in the
  obvious way.
\end{proof}

\begin{comment}
\begin{lemma} Let $\coprod:\sC\times\sC\to\sC$ be a bifunctor provided
with the following data:
\begin{itemize}
\item a symmetric monoidal structure with unit object $e$,
\item two ``coprojections" $c_1^{A,B}:A\to A\coprod B$, $c_2^{A,B}:B\to
A\coprod B$,
natural in $A$ and $B$,
\item a ``unit" $\eta^A:e\to A$, natural in $A$,
\item a ``codiagonal'' $\nabla^A:A\coprod A\to A$, natural in $A$,
\end{itemize}
satisfying the following identities:
\begin{gather*}
\nabla_A\circ c_1^{A,A} =\nabla_A\circ c_2^{A,A} = 1_A\\
\nabla_{A\coprod B}\circ \nabla_{A\coprod B\coprod A\coprod B}\circ
(c_1^{A,B\coprod A\coprod B}\coprod c_2^{A\coprod B\coprod A,B}) =
1_{A\coprod B},
\end{gather*} 
the second one modulo the associativity constraint. Then $\coprod$ is a
binary coproduct in $\sC$. 
\end{lemma}

\begin{proof} Given two morphisms $f:A\to C$, $g:B\to C$, define
$(f,g):A\coprod B\to C$ by
\[(f,g) = \nabla_C\circ (f\coprod g).\]

The first identity implies that $(f,g)\circ c_1^{A,B} = f$ and $(f,g)\circ
c_2^{A,B} = g$. It remains to show that if $h:A\coprod B\to C$ verifies
the same properties, then $h = (f,g)$. But, from the second identity, we
get
\begin{multline*}
h=h\circ \nabla_{A\coprod B}\circ \nabla_{A\coprod B\coprod A\coprod B}\circ
(c_1^{A,B\coprod A\coprod B}\coprod c_2^{A\coprod B\coprod A,B})\\
= \nabla_C\circ (h\coprod h)\circ \nabla_{A\coprod B\coprod A\coprod
B}\circ (c_1^{A,B\coprod A\coprod B}\coprod c_2^{A\coprod B\coprod A,B})
\end{multline*}
\end{proof}
\end{comment}

\begin{proposition}[\protect{\cite[1.3.6 and 2.1.8]{maltsin}}]\label{p6.2}
  Let $\sC$ be a 
  category with finite coproducts and $S$ a 
  family of morphisms of 
  $\sC$ stable under coproducts.
Then $S^{-1}\sC$ has finite
  coproducts and the localisation functor $\sC\to S^{-1}\sC$ commutes
  with them.\qed
\end{proposition}

\begin{corollary} Let $\sC$ be a category and $S$ a
  family of morphisms of   $\sC$. In $\sC^\Coprod$, consider the
  following family $S^\Coprod$ (see proof of Proposition \ref{p6.1}): 
  $s:(C_i)_{i\in I}\to (D_j)_{j\in J}$ is in $S^\Coprod$ if and only
  if the underlying map $f:I\to J$ is bijective and $s_i:C_i\to
  D_{f(i)}$ belongs to $S$ for all $i\in I$.
Then we have an equivalence of categories 
\[(S^{-1}\sC)^\Coprod\simeq(S^\Coprod)^{-1}\sC^\Coprod.\] 
\end{corollary}

\begin{proof} By Proposition \ref{p6.2},
  $(S^\Coprod)^{-1}\sC^\Coprod$ has finite coproducts, hence it is
  enough to show that any functor 
  $F:S^{-1}\sC\to\sE$, where $\sE$ has finite coproducts, factors
  canonically through a functor from $(S^\Coprod)^{-1}\sC^\Coprod$
  which commutes with finite coproducts. Let $P:\sC\to
  S^{-1}\sC$ be the localisation functor; then $F\circ P$ factors
  through $\sC^\Coprod$. The resulting functor inverts morphisms
  of $S$ and commutes with coproducts, hence also inverts morphisms of
  $S^\Coprod$. Thus we get a functor
  $(S^\Coprod)^{-1}\sC^\Coprod\to\sE$, which obviously commutes with
  finite coproducts. 
\end{proof}

\begin{theorem} \label{t6.1} In the situation of \eqref{eq0}, if $\bar T$ is an
  equivalence of categories, then so is $\overline{T^\Coprod}$, where
  $T^\Coprod:\sC^\Coprod\to \sD^\Coprod$ is the functor induced by
  $T$. Moreover,  $\overline{T^\Coprod}=(\bar T)^\Coprod$.\qed
\end{theorem}

\section{Applications in algebraic geometry}\label{7}

Let $k$ be a field. We denote by $\Sch(k)$ the category of reduced separated
$k$-schemes of finite type. 

\subsection{Hyperenvelopes (Gillet--Soul\'e \cite{gs})} In this
example, $k$ is of 
characteristic
$0$. We take for $\sD^\op$ the category of simplicial reduced
$k$-schemes of finite 
type, and for $\sC^\op$ the full subcategory consisting of smooth simplicial
$k$-schemes.

For $S$ and $S'$ we take \emph{hyperenvelopes} as considered by Gillet
and Soul\'e 
in \cite[1.4.1]{gs}: recall that a map $f:X_\bullet\to Y_\bullet$ in $\sD$ is a
hyperenvelope if and only if, for any extension $F/k$, the induced map
of simplical 
sets $X_\bullet(F)\to Y_\bullet(F)$ is a trivial Kan fibration (see
loc. cit. for 
another equivalent condition).

\begin{theorem}\label{th} In the above situation, the conditions of
Theorem \ref{p3} are satisfied. In particular, $\bar T^E$ is an equivalence of
categories for any finite ordered set $E$.
\end{theorem}

\begin{proof} (i) is true by definition; (ii) is proved (or remarked) in
\cite[p. 136]{gs} and (iv) is proved in \cite[Lemma 2 p. 135]{gs}
(which, of course, 
uses Hironaka's resolution of singularities). The last assertion
follows from Theorem 
\ref{p3}.
\end{proof}

\subsection{Proper hypercovers (Deligne--Saint Donat
  \protect\cite{SGA4.2})} Here $k$ 
is any field. We take the same
$\sC$ and
$\sD$ as in the previous example, but we let $S'$ be the collection of proper
hypercovers (defined from proper surjective morphisms as in \cite[Exp.
Vbis, (4.3)]{SGA4.2}). 

\begin{theorem}\label{tp} In the above situation, the conditions of
Theorem \ref{p3} are satisfied. In particular, $\bar T^E$ is an equivalence of
categories for any finite ordered set $E$.
\end{theorem}

The proof is exactly the same as for Theorem \ref{th}, replacing the
use of Hironaka's theorem in the proof of (iv) by that of de Jong's
alteration theorem 
\cite{dJ}.

\subsection{Cubical hyperresolutions (Guill\'en--Navarro Aznar
  \protect\cite{GNPP})} In 
this example,
$k$ is again of characteristic $0$. We are not going to give a new proof of
the main theorem 
of \cite[Th. 3.8]{GNPP}, but merely remark that its proof in
loc. cit. can be viewed as checking a special case of Theorem 
\ref{t0}. Namely, in this situation, $\sD$ is a category of diagrams
of schemes of a certain type, $\sC$ is the category of cubical
hyperresolutions of objects of $\sD$, $T$ maps a hyperresolution to
the diagram it resolves, $S'$ consists of identities and $S$ consists
of 
arrows mapping to identities; the categories $I_d$, $I_f$ then reduce
to the fibre 
categories of $T$. Guill\'en and Navarro Aznar prove that, on the
level of $S^{-1}\sC$, $I_d$ is 
$1$-connected for any $d\in\sD$ and that $I_f$ is $0$-connected for any $f\in
\sD^{\Delta^1}$. The $-1$-connectedness of the $I_{f_2,f_1}$ is then
automatic in this 
special case, because Lemma 3.8.6 of \cite{GNPP} shows that the first
two conditions 
already imply that $\bar T$ is faithful.

\subsection{Jouanolou's device (Riou \protect\cite[Prop. II.16]{riou})}
Here $\sC$ is the category of smooth affine schemes over some
regular scheme $R$, $\sD$ is the category of smooth $R$-schemes,
$S'$ consists of morphisms of the form $Y\to X$ where $Y$ is a torsor
under a vector bundle on $X$ and $S = S'\cap \sC$. Riou checks that
the hypotheses of Theorem \ref{riou} are verified by taking the opposite
categories, hence that the inclusion functor $T:\sC\to\sD$ induces an
equivalence on localised categories.

\subsection{Closed pairs} Here we take for $\sC$ the category whose
objects are closed 
embeddings $i:Z\to X$ of proper $k$-schemes such that $X-Z$ is dense
in $X$, and a 
morphism from
$(X,Z)$ to
$(X',Z')$ is a morphism $f:X\to X'$ such that $f(X-Z)\subseteq X'-Z'$. We take
$\sD=\Sch(k)$, and for
$T$ the functor
$T(X,Z)=X-Z$. Finally, we take for $S'$ the isomorphisms of $\sD$ and
$S:=T^{-1}(S')$. 

\begin{theorem}[{\it cf.} \protect{\cite[Lemma 2.3.4]{GN2}}]\label{t10} In the above 
situation, the conditions  of
Proposition \ref{p1} b) are satisfied. In particular, $\bar T$ is an
equivalence of 
categories.
\end{theorem}

\begin{proof} It is sufficient to check (b1) and the
conditions of Proposition \ref{p1} c). In (b1),
$T(f)s=T(g)s\Rightarrow T(f)=T(g)$ is trivial since $s$ is by definition an
isomorphism. On the other hand, $T$ is faithful by a classical
diagonal argument, since 
all schemes are separated.

In Proposition \ref{p1} c), the assertion on finite products is clear
(note that 
$(X_1,Z_1)\times (X_2,Z_2)= (X_1\times X_2,Z_1\times X_2\cup X_1\times
Z_2)$). For 
$U\in \Sch(k)=\sD$, we define $K_U$ as the full subcategory of $J_U$
consisting of 
immersions $U\inj X-Z$.

Nagata's theorem implies that $I_U$ is nonempty; in particular, $K_U$
is nonempty. Let 
$\kappa=(U\inj X-Z)$ be an object of $K_U$, and let $\bar U$ be the
closure of $U$ in $X$. 
Then $(\bar U,\bar U-U)$ defines an object of $I_U/\kappa$, and (c1)
is verified. As for 
(c2), it is trivial since the product of an immersion with any
morphism remains an 
immersion.
\end{proof}

\subsection{Another kind of closed pairs} Here we assume that $\car k
= 0$. For $n\ge 
0$, we define
$\sD_n^\op$ to be the category whose objects are closed embeddings
$i:Z\to X$ with $X$ an (irreducible) variety of dimension $n$, 
$X-Z$ dense and
\emph{smooth}; a morphism from $(X,Z)$ to $(X',Z')$ is a  map $f:X\to
X'$ such that $f^{-1}(Z')=Z$.  
We define
$\sC_n^\op$ as the full subcategory of
$\sD_n^\op$ consisting of pairs $(X,Z)$ such that $X$ is smooth.

We take for $S'$ the set of morphisms $s:(X,Z)\to (X',Z')$ such that
$s_{|X-Z}$ is an 
isomorphism onto $X'-Z'$, and $S=S'\cap \sC_n$.

\begin{lemma}\label{l1} If $f$ and $s$ have the same domain in
  $\sD_n$, with $s\in S'$, then $s$ is always in good position with
  respect to $f$.  
\end{lemma} 

\begin{proof} Translating in the opposite category, we have to see
  that if $f:(X_1,Z_1)\to (X,Z)$ and $s:(\tilde X,\tilde Z)\to (X,Z)$
  are maps in $\sD^\op$ with $s\in S'$, then the fibre product
  $(\tilde X_1,\tilde Z_1)$ of $f$ and $s$ exists and the pull-back
  map $s':(\tilde X_1,\tilde Z_1)\to (X_1,Z_1)$ is in $S'$. Indeed,
  note that $(\tilde X_1,\tilde Z_1)$ is given by the same formula as
  in the proof of Theorem \ref{t10} provided it exists, namely $\tilde
  X_1 = X_1\times_X \tilde X$ and $\tilde Z_1 = X_1\times_X \tilde
  Z\cup Z_1\times_X \tilde X$. The thing to 
check is that $\tilde X_1 - \tilde Z_1$ is still dense in $\tilde
  X_1$, which will 
imply in particular that $\tilde X_1$ is a variety. It is
sufficient to check separately that $\tilde Z\times_X X_1$ and $\tilde
  X\times_X 
Z_1$ are nowhere dense in $\tilde X_1$, which we leave to the reader.
\end{proof}

\begin{theorem}\label{t8} In the above situation, the conditions of
  Proposition \ref{p2} are 
satisfied. In particular, $\bar T$ is an equivalence of
categories. Moreover, we don't change $S^{-1}\sC$ if we replace $S$ by
the subset of $S'\cap \sC_n$ generated by blow-ups with smooth
centres. 
\end{theorem}

\begin{proof} (d1) is true because two morphisms from the same source
  which coincide on a dense open subset are equal. (d2) and (d3) are
  immediately checked thanks to Lemma \ref{l1}. (d4) is clear and (d5)
  follows from Hironaka's resolution theorem. The last statement of
  Theorem \ref{t8} also follows from Hironaka's theorem that any
  resolution of singularities may be dominated by a composition of
  blow-ups with smooth centres. 
\end{proof}

\begin{comment}
\subsection{Factoring morphisms} This has applications in duality 
constructions. Let $X$ be the base scheme and $\Sch(X)$ denote the 
category of $X$-schemes of finite type. We denote by $\Sch^2_{cs}(X)$ 
the 
category whose objects are closed immersions $i:Y\to Z$ in $\Sch(X)$, with $f:Z\to X$ smooth, and morphisms are commutative squares. Similarly, we denote by (resp. $\Sch^2_{op}(X)$) the 
category whose objects are open immersions $j:Y\to Z$ in $\Sch(X)$, with $p:Z\to X$ proper, and morphisms are commutative squares. We have forgetful 
functors
\[T_{cs}:\Sch^2_{cs}(X)\to \Sch(X)\]
\[T_{op}:\Sch^2_{op}(X)\to \Sch(X)\]
sending $(Y\to Z)$ to $Y$.

Let $S'\subset Ar(\Sch(X))$ denote the identity morphisms and $S_{cs}=T_{cs}^{-1}(S')$, $S_{op}=T_{op}^{-1}(S')$. Then it is easily checked 
that the conditions of Theorem \ref{riou} hold and hence Theorem 
\ref{p3} along with Lemma \ref{inftycon} can be applied to obtaint he 
following result. 
\begin{theorem}
The categories $S^{-1}\Sch^{cs}(X)$ and $S^{-1}\Sch^{op}(X)$ are 
equivalent to $\Sch(X)$.
\end{theorem}
\end{comment}

\section{Applications in birational geometry}\label{8}

We shall reserve the word ``variety" to mean an integral
scheme in $\Sch(k)$, and denote their full subcategory by $\Var(k)$; we
usually abbreviate with $\Sch$ and $\Var$. Recall 
\cite[(2.3.4)]{ega} that a
\emph{birational morphism}
$s:X\to Y$ in $\Sch$ is a morphism such that every irreducible component $Z'$
of $Y$ is dominated by a unique irreducible component $Z$ of $X$ and
the induced map 
$s_{|Z}:Z\to Z'$ is a birational map of varieties. 

\begin{definition}\label{d1}
We denote by $S_b$ the multiplicative system of birational morphisms in
$\Sch$, by $S_o$ the subsystem consisting of open immersions and by $S_b^p$ the
subsystem consisting of proper birational morphisms.
\end{definition}

We shall also say that a morphism $f:X\to Y$ in $\Sch$ is \emph{dominant} if its
image is dense in $Y$, or equivalently if every irreducible component of $Y$ is
dominated by some irreducible component of $X$.

\begin{lemma}\label{l5} a) Let 
\[\begin{CD}
X&\begin{smallmatrix}\sigma_1\\
  \rightrightarrows\\\sigma_2\end{smallmatrix}& Y\\ 
&& @V{s}VV\\
&& Z.
\end{CD}\]
be a diagram in $\Sch$, with $X$ reduced, $Y$ separated and
$s,\sigma_1,\sigma_2\in S_b$. Suppose that
$s\sigma_1=s\sigma_2$. Then $\sigma_1=\sigma_2$.\\
b) Let  $f,g:Y\to Z$, $h:X\to Y\in Ar(\Sch)$ be such that $fh=gh$. Suppose that
$h$ is dominant. Then $f=g$.
\end{lemma}

\begin{proof} a) Recall from \cite[(5.1.5)]{ega} the kernel scheme
$\Ker(\sigma_1,\sigma_2)\subset 
X$: it is the  inverse image scheme of the diagonal
$\Delta_Y(Y)\subset Y\times_k Y$ 
via the morphism $(\sigma_1,\sigma_2)$. Since $Y$ is separated,
$\Ker(\sigma_1,\sigma_2)$ is a closed subscheme of $X$ and, by definition of
birational morphisms, it contains all the generic points of $X$. Hence
$\Ker(\sigma_1,\sigma_2)=X$ since $X$ is reduced, and $\sigma_1=\sigma_2$.

b) is obvious, since by assumption $h^{-1}(\Ker(f,g))=X$.
\end{proof}

\begin{definition} Let $\sC$ be a subcategory of $\Sch$.\\
a) We denote by $\sC^\qp$ (resp. $\sC^\proper,\sC^\proj$) the full
subcategory of $\sC$ consisting of quasiprojective (resp. proper,
projective) objects.\\ 
b) We denote by $\sC_\sm$ the non-full subcategory of $\sC$ with the
same objects, but where a morphism $f:X\to Y$ is in $\sC_\sm$ if and
only if $f$ maps the smooth locus of $X$ into the smooth locus of
$Y$. 
\end{definition}

The following proposition is the prototype of our birational results.

\begin{proposition} \label{p3.5} In the commutative diagram
\[\xymatrix{
S_b^{-1}\Var^\proper\ar[r]^A& S_b^{-1}\Var\\
S_b^{-1}\Var^\proj\ar[r]^B\ar[u]_C& S_b^{-1}\Var^\qp\ar[u]_D
}\]
all functors are equivalences of categories. The same holds by adding the
subscript $\sm$ everywhere.
\end{proposition}

\begin{proof} We first prove that $A$ and $B$ are equivalences of
categories. For this, we apply Proposition \ref{p1} b) with
$\sC=\Var^\proper$ (\resp $\Var^\proj$),
$\sD=\Var$ (\resp $\Var^\qp$), $T$ the obvious inclusion,
$S=S_b$ and $S'=S_b$:
\begin{itemize} 
\item Condition (b1) holds because $T$ is fully
  faithful and birational morphisms are dominant (see Lemma \ref{l5} b)).
\item (b2) is true by Nagata's Theorem in the proper case and
  tautologically in the projective case.
\item For (b3) we use the ``graph trick'': we are given $i:X\to\bar X$
  and $j:X\to Y$ where $\bar X$ and $Y$ are proper (resp. projective)
  and $i$ is 
  birational. Let $\bar X'$ be the closure of the diagonal image of
  $X$ in $\bar X\times Y$: then $X\to\bar X'$ is still birational,
  $\bar X'$ is proper (resp. projective) and the projections $\bar
  X'\to\bar X$, $\bar X'\to Y$ give the desired object of
  $I_X/i\times_{I_X} I_X/j$.
\end{itemize}

We now prove that $D$ is an equivalence of categories, which will also
imply that $C$ is an equivalence of categories. This time we apply
Proposition \ref{p2} with
  $\sC=(\Var^\qp)^\op$, $\sD=\Var^\op$, $T$ the
  obvious inclusion and $S=S_o$, $S'=S_o$: 
\begin{itemize}
\item Condition (d1) is clear (open immersions are monomorphisms even
  in $\Var$). 
\item (d2) means that the intersection of two dense open subsets in a
  variety is dense, which is true. 
\item (d3) means that if $(gf)^{-1}(U)\neq\emptyset$, then
  $g^{-1}(U)\neq\emptyset$, which is true. 
\item (d4) is clear.
\item In (d5), we have a morphism $f:X_1\to X$ of varieties and want to find a
quasi-projective dense open subset $U\subseteq X$ such that
$f^{-1}(U)\neq \emptyset$: take $U$ containing $f(\eta_{X_1})$ (any
point has an affine neighbourhood). 
\end{itemize}

The proofs with indices $\sm$ are the same.
\end{proof}

\begin{proposition} \label{p3.6} In the commutative diagram
\[\xymatrix{
S_b^{-1}\Sm^\proper\ar[r]^{A'}& S_b^{-1}\Sm\\
S_b^{-1}\Sm^\proj\ar[r]^{B'}\ar[u]_{C'}& S_b^{-1}\Sm^\qp\ar[u]_{D'}
}\]
$D'$ is an equivalence of categories. Under resolution of
singularities, this is true of the three other functors.  
\end{proposition}

\begin{proof} The same as that of Proposition \ref{p3.5}, except that
  for $A'$ and 
$B'$, we need to desingularise a compactification of a smooth variety
using Hironaka's Theorem.
\end{proof}

\begin{proposition} \label{p3.7} If $k$ is perfect, in the commutative diagram
\[\xymatrix{
S_b^{-1}\Sm\ar[r]^E & S_b^{-1}\Var_\sm\\
S_b^{-1}\Sm^\qp\ar[r]^F \ar[u]_{G} &S_b^{-1}\Var_\sm^\qp\ar[u]_H
}\]
all functors are equivalences of categories.
\end{proposition}

\begin{proof} The case of $H$ has been seen in Proposition \ref{p3.5}, and the
case of $G=D'$ has been seen in Proposition \ref{p3.6}.  We now prove
that $E$ and
$F$ are equivalences of categories. Here we apply Proposition \ref{p2} with
$\sC=\Sm^\op$ (\resp $(\Sm^\qp)^\op$),  $\sD=
\Var_\sm^\op$ (\resp $(\Var_\sm^\qp)^\op$), $T$ the obvious
inclusion and $S=S'=S_o$. Note that open immersions automatically
respect smooth loci. Let us check the conditions: 
\begin{itemize}
\item (d1), (d2) and (d3) and (d4) are clear (see proof of Proposition
  \ref{p3.5}). 
\item It remains to check (d5): if $f:X_1\to X$ is a morphism in
  $\Var_\sm$, then $f(\eta_{X_1})$ is contained in the smooth locus
  $U$ of $X$, hence $U\to X$ is in good position with respect to $f$.  
\end{itemize}
\end{proof}

\begin{proposition} \label{p3.8} Under resolution of singularities, all
  functors in the
  commutative diagram
\[\xymatrix{
S_b^{-1}\Sm^\proper\ar[r]^I&
S_b^{-1}\Var_\sm^\proper\\ 
S_b^{-1}\Sm^\proj\ar[r]^J\ar[u]_K& S_b^{-1}\Var_\sm^\proj\ar[u]_L
}\]
are equivalences of categories.
\end{proposition}

\begin{proof} The case of $K=C'$ has been seen in Proposition
  \ref{p3.6} and the 
case of $L$ in Proposition \ref{p3.5}. The case of the other functors is
then implied by the previous propositions (the reader should draw a
commutative cube of categories in order to check that enough equivalences
of categories have been proven).
\end{proof}

\begin{proposition} \label{p8.1} The previous propositions remain true if we
  replace all 
  categories in sight by their finite coproduct envelopes (see
  Proposition \ref{p6.1}) and $S_b$ by $S_b^\Coprod$ (ibid.).
\end{proposition}

\begin{proof} This follows from Theorem \ref{t6.1}.
\end{proof}

\begin{remarks} a) Note that even though Proposition \ref{p8.1} says that
  $D^\Coprod$ induces an equivalence of categories on localisations,
  where $D$ is the functor of Proposition \ref{p3.5},
  $(D^\Coprod,S_o^\Coprod)$ does 
  not satisfy the (dual) simplicial hypotheses of Theorem \ref{t0}. Indeed,
  let $X$ be a non-quasiprojective variety over $k$ that we assume
  algebraically closed for simplicity. By Kleiman's theorem
  \cite{kleiman}, there exists a finite set $\{x_1,\dots,x_n\}$ of
  closed points 
  of $X$ which is contained in no affine open
  subset, hence also in no quasi-projective open subset. Thus, if
  $Y=\coprod_n \Spec k$ and $f:Y\to X$ is the map defined by the
  $x_i$, then $I_f$ is empty. This shows that the simplicial
  hypotheses are not preserved by finite product envelope.\\
b) Also, while $(D,S_o)$ satisfies the dual simplicial hypotheses, it
does not satisfy the dual of hypothesis (1') of Corollary \ref{c2}:
this is obvious from Chow's lemma. This shows that the hypotheses of
Corollary \ref{c2} are strictly stronger than those of Theorem \ref{t0}. 
\end{remarks}

\begin{remark}\label{r3.2} To summarise Propositions \ref{p3.6},
  \ref{p3.7} and \ref{p3.8} under resolution of singularities, we have
  the following equivalences of categories: 
\begin{multline*}
S_b^{-1}\Sm^\proj\simeq S_b^{-1}\Sm^\proper\simeq S_b^{-1}\Sm^\qp
\simeq S_b^{-1}\Sm \\ 
\simeq S_b^{-1}\Var_\sm^\proj \simeq S_b^{-1}\Var_\sm^\proper\simeq
S_b^{-1}\Var_\sm^\qp\simeq S_b^{-1}\Var_\sm. 
\end{multline*}

(One could also replace the superscript $\qp$ by ``affine", as the proofs
show.)

We shall show in \cite{ks} that
\[S_b^{-1}\Sm^\proj(X,Y) = Y(k(X))/R\]
for any two smooth projective varieties $X,Y$,  where $R$ is Manin's
$R$ equivalence.
\end{remark}

\begin{remark} On the other hand, the functor
$S_b^{-1}\Sm\to
  S_b^{-1}\Var$ is neither full nor faithful, even under resolution of
singularities. Indeed, take $k$ of characteristic $0$. Let $X$ be a
  proper irreducible curve of geometric
  genus $>0$ with one nodal singular point $p$. Let $\bar \pi:\tilde
  X\to X$ be its normalisation, $U=X-\{p\}$, $\tilde U=\bar
  \pi^{-1}(U)$, $\pi=\bar \pi_{|\tilde U}$ and $j:U\to X$, $\tilde
\jmath:\tilde
  U\to \tilde X$ the two inclusions. We assume that $\bar \pi^{-1}(p)$
  consists of two rational points $p_1,p_2$. Finally, let $f_i:\Spec
  k\to \tilde X$ be the map given by $p_i$. 
\[\xymatrix{
\Spec F\ar@< 2pt>[r]^{f_1}\ar@<-2pt>[r]_{f_2}& \tilde X\ar[r]^{\bar
\pi}&X\\ &\tilde U\ar[u]^{\tilde \jmath}\ar[r]^\pi& U\ar[u]_j
}\]

In $S_b^{-1}\Var$,  $\bar \pi$ is an isomorphism so that
$f_1=f_2$. We claim that $f_1\ne f_2$ in
$S_b^{-1}\Sm^\proper\iso S_b^{-1}\Sm$. Otherwise, since
$R$-equivalence is a birational invariant of smooth proper
varieties \cite[Prop. 10]{ct-s}, we would have $p_1=p_2\in \tilde
X(k)/R$. But this is false
since $\tilde X$ has nonzero genus. We thank A. Chambert-Loir for his help in
finding this example.

More generally, it is well-known that for any integral curve $C$ and
any two closed points
$x,y\in C$ there exists a proper birational morphism $s:C\to C'$
such that
$s(x)=s(y)$ (\cf \cite[Ch. IV, \S 1, no 3]{serre2} when $F$ is
algebraically closed). This shows that any two morphisms
$f,g:X\rr C$ such that
$f(\eta_X)$ and $g(\eta_X)$ are closed points become equal in
$S_b^{-1}\Var$. This can be used to show that the functor
$S_b^{-1}\Sm\to S_b^{-1}\Var$ \emph{does not have a right adjoint}.

Non fullness holds even if we restrict to normal varieties. Indeed,
let us take $k=\R$ and let $X$ be the affine cone with
equation $x_1^2+\dots +x_n^2=0$ (for $n\ge 3$ this is a normal
variety). Let $\tilde X$ be a desingularisation of $X$ 
(for example obtained by blowing up the singular point) and  $\bar X$
a smooth compactification of $\tilde X$. Then $\bar X(\R)=\emptyset$
by a valuation argument, 
hence $S_b^{-1}\Sm^\proj(\R)(\Spec \R,\bar X)=\emptyset$ by Remark
\ref{r3.2}. On the other
hand, $X(\R)\neq\emptyset$, hence 
\[S_b^{-1}\Var^\proj(\R)(\Spec
\R,\bar X)=S_b^{-1}\Var^\proj(\R)(\Spec \R,X)\neq\emptyset.\] 
We are
indebted to Mah\'e for pointing out this example. For $n\ge 4$, this
singularity is 
even terminal in the sense of Mori's minimal model programme, as
Beauville pointed out (which seems to mean unfortunately that we
cannot insert this programme in our framework...)  
\end{remark}

\begin{remark} Let $n\ge 0$. Replacing all the subcategories $\sC$ of
  $\Sch$ used above 
  by their full subcategories $\sC_n$ consisting of schemes of
  dimension $\le n$, one checks readily that all corresponding
  equivalences of categories
  remain valid, with the same proofs. This raises the question whether
  the induced functor $S_b^{-1}\sC_n\to S_b^{-1}\sC_{n+1}$ is fully
  faithful for some (or all) choices of $\sC$. It can be proven
  \cite{ks} that this is true at least for $\sC=\Sm^\proj$ in
  characteristic zero, hence for the other $\sC$s which become
  equivalent to it after inverting birational morphisms as in Remark
  \ref{r3.2}. However, the
  proof is indirect and consists in observing that the morphisms
  are still given by the formula of Remark \ref{r3.2}. It is an
  interesting question
  whether such a result can be proven by methods in the spirit of the
  present paper. 
\end{remark}

\end{document}